\documentclass{amsart}
\usepackage{amsmath}
\usepackage{amssymb}
\usepackage{euscript}
\usepackage[all]{xy}
\usepackage{varioref}
\usepackage{verbatim}
\usepackage{xspace}
\newtheorem{theorem}{Theorem}
\newtheorem{corollary}{Corollary}
\newtheorem{lemma}{Lemma}

\newtheorem{definition}{Definition}

\newcommand{\pd}[2]{\ensuremath{\frac{\partial{#1}}{\partial{#2}}}}
\newcommand{\R}[1]{\ensuremath{\mathbb{R}^{#1}}}
\newcommand{\C}[1]{\ensuremath{\mathbb{C}^{#1}}}

\newcommand{\CP}[1]{\ensuremath{\mathbb{CP}^{#1}}}

\newcommand{\Gro}[2]{\widetilde{\operatorname{Gr}}\left({#1},{#2}\right)}

\newcommand{\GL}[1]{\operatorname{GL}\left({#1}\right)}

\begin{document}
\title{Almost complex rigidity of the complex projective plane}
\author{Benjamin McKay}
\address{University College Cork \\
Cork, Ireland}
\email{B.McKay@UCC.ie}
\date{\today \\ MSC Primary: 57R17, 51A10; Secondary: 14N05}
\begin{abstract} An isomorphism of
symplectically tame smooth pseudocomplex structures 
on the complex projective
plane which is a homeomorphism and 
differentiable of full rank at two points is
smooth.
\end{abstract}
\thanks{This work was supported in full or in part by a grant from 
the University of South Florida St. Petersburg New Investigator
Research Grant Fund. This support does not necessarily imply
endorsement by the University of research conclusions.}
\maketitle
\tableofcontents
\section{Introduction}
A longstanding concern of Gromov and of Uhlenbeck 
(among many other mathematicians) is to understand global features
of elliptic partial differential equations using global geometry and 
weak local estimates.
This paper employs global symplectic geometry and 
weak local estimates 
to prove smoothness of isomorphisms of pseudocomplex
structures on the complex projective plane. 

A \emph{pseudocomplex structure} on a 4 dimensional
manifold $M$ is a choice, for any smooth complex-valued 
local coordinate system $z,w : M \to \C{}$, 
of a system of partial differential equations
\[
\pd{w}{\bar{z}} = 
F\left(z,\bar{z},w,\bar{w},\pd{w}{z},\pd{\bar{w}}{\bar{z}}\right),
\]
(with $F$ a smooth function), so that under smooth changes of 
coordinates, the partial differential equations are equivalent,
i.e. have the same local solutions. 
An example: every almost
complex structure has Cauchy-Riemann equations for pseudoholomorphic
curves, giving a pseudocomplex structure. 
If $E$ is the pseudocomplex
structure, the local solutions of the partial differential 
equations are smooth surfaces in the manifold $M$, called
\emph{$E$-curves}. The concept is due to Gromov \cite{Gromov:1985},
pg.~342.
McKay \cite{McKay:1999,McKay:2003} analyzed
$E$-curves locally and globally. 

A \emph{morphism} of
pseudocomplex structures $E_0$ on $M_0$ and $E_1$ on $M_1$
is a map $M_0 \to M_1$ carrying $E_0$-curves to $E_1$-curves.
A pseudocomplex structure $E$ on a 4 dimensional manifold $M$ 
is \emph{tamed} by a symplectic
structure $\omega$ on $M$ if $\omega>0$ on all $E$-curves,
it is \emph{tame} if it is tamed by some symplectic structure.
\begin{theorem}[Main theorem]
A isomorphism of smooth tame pseudocomplex structures
on the complex projective plane which is a homeomorphism, and
differentiable of full rank at two points,
is smooth.
\end{theorem} 
\begin{corollary} A continuous biholomorphic map between
smooth tame almost complex structures on the complex 
projective plane which
is differentiable of full rank at two points, is smooth.
\end{corollary}
\begin{corollary}
A biLipschitz isomorphism of smooth tame pseudocomplex
structures on the complex projective plane is smooth.
\end{corollary}
The proof is purely geometric, using
elementary projective plane geometry coupled with arguments from
Gromov \cite{Gromov:1985}. By comparison, the strongest known
result for almost complex 4 dimensional manifolds is due to
Coupet, Gaussier \& Sukhov \cite{CGS:2003}, proving smoothness
assuming continuous differentiability (although their result 
extends to a pseudoconvex boundary).
\section{Generalities on pseudocomplex structures}
See McKay \cite{McKay:1999,McKay:2003} for an extensive discussion of
pseudocomplex structures. 
Let us give a slightly more geometric definition of pseudocomplex
structures, following Gromov. 
\begin{definition}
A \emph{pseudocomplex structure} on a four dimensional manifold $M$ is a
choice of smooth immersed submanifold $E \subset \Gro{2}{TM}$ 
inside the bundle of oriented 2-planes in the tangent spaces 
of $M$, so that the map $E
\to M$ is a submersion, and so that the requirement that a surface $C
\subset M$ have tangent planes belonging to $E$ be equivalent in local
coordinates to a determined elliptic system of partial differential 
equations.
\end{definition} 
As Gromov points out, and my previously cited articles prove, the
ellipticity requirement can be expressed neatly in terms of the
canonical $(2,2)$ signature conformal structure of the Grassmannians 
$\Gro{2}{T_m M}=\Gro{2}{\R{4}}$ (invariant under $\GL{4,\R{}}$) as the
requirement that the fibers $E_m \subset \Gro{2}{T_mM}$ 
of $E$ are definite surfaces in that
conformal structure, i.e. nowhere tangent to null directions.
We will refer to a manifold $M$ with pseudocomplex structure $E$ as a
pseudocomplex manifold.
\begin{definition}
An oriented immersed surface $C \subset M$ in a pseudocomplex 
manifold is called an \emph{$E$-curve} if its tangent spaces 
belong to $E$.
\end{definition}
There are two natural notions of morphism of a pseudocomplex
structure: we could require, as we have chosen to do, only that a
morphism $\left(M_0,E_0\right) \to \left(M_1,E_1\right)$ 
take $E_0$-curves to $E_1$-curves, or we could require that it also
preserve orientations. We will ignore the orientations. For example,
the automorphism group of the usual complex structure on $\CP{2}$
includes conjugate holomorphic maps, given in local affine charts by
complex conjugation.
\begin{theorem}[McKay \cite{McKay:1999,McKay:2003}]
Let $E \subset \Gro{2}{TM}$ be a pseudocomplex structure.
Pick a point $e \in E$ and let $m \in M$ be its image under 
$E \to \Gro{2}{TM} \to M.$
To each point $e \in E$ there is associated smoothly a choice of 
\emph{osculating complex structure} $J_e : T_m M \to T_m M.$ 
The 2-plane $e$ is a $J_e$-complex line.
In particular, if $C$ is an $E$-curve, then at each point $m \in C$,
$e=T_m C$ is invariant under $J_e$. In particular, $C$ inherits the
structure of a Riemann surface. The osculating complex structure
is uniquely determined by the requirement that the linearization of
the elliptic equation for $E$-curves at any $E$-curve $C$ is 
of the form
\[
\bar{\partial} \sigma = \tau\left( \sigma \right )
\]
where $\sigma$ is a section of the normal bundle of $C$, 
$\tau$ is a $J_e$-conjugate linear map, and $J_e$ is used
to fix the meaning of $\bar{\partial}.$ (This $\tau$ is 
construced from the torsion functions 
$T_2, T_3$ of McKay \cite{McKay:2003}. We will refer to
a solution $\sigma$ of this equation as a 
\emph{pseudoholomorphic}
normal vector field.)
This linearized equation on sections of the normal bundle 
makes the $E$-curve $C$ into a 
Riemann surface, and the normal
bundle into a holomorphic line bundle.
\end{theorem}
\begin{proof}
Most of this is proven in McKay \cite{McKay:2003}.
The remarks on the linearization can be easily checked
from the structure equations in that paper:
a one parameter family of $E$-curves, say with parameter 
$t$,  will lift to $E$ to
look like (in the notation and terminology of 
that paper) $\theta = s \, dt.$ This $s$ function on the
bundle of adapted frames over the total space of the 
family of $E$-curves varies under the representation of
the structure group (as the reader can check) which is
the defining representation of the normal bundle.
We linear just by taking exterior derivative, and then
setting $t=0$. We can always adapt frames (as in that
paper) to arrange $\pi=0,$ so that we find by taking
the exterior derivative: $\bar{\partial}s + \tau_1 \bar{s}=0.$
Check that this equation descends from the bundle
of adapted frames down to the $E$-curve itself. 
The complex structure of the $E$-curve is 
determined by the characteristic variety of this
equation---the characteristic variety is precisely
the union of the holomorphic and conjugate holomorphic
tangent bundles. The orientation of the $E$-curve
picks out the holomorphic tangent bundle. 
By following Duistermaat \cite{Duistermaat:1972}, 
we find the complex structure $J_e$ completely
determined, and indeed the structure of holomorphic
line bundle on the normal bundle, from the form of 
the differential equation above for pseudoholomorphic
sections of the normal bundle.  Note that the
torsion prevents the first-order deformations
of the $E$-curve being identified with the holomorphic
sections of this line bundle.
Nonetheless, as
Duistermaat points out, we can carry out Riemann--Roch/Chern class
theory of these $\bar{\partial} \sigma + \tau\left(\sigma\right)=0$
equations just as for Cauchy--Riemann equations, since
the behaviour near zeroes of $\sigma$ is identical.
\end{proof} 
\begin{definition}
A pseudocomplex structure $E \subset \Gro{2}{TM}$ 
is \emph{proper} if its fibers $E_m \subset \Gro{2}{T_mM}$
are compact (in other words, if the map $E \to M$ is a proper map).
\end{definition}
All pseudocomplex structures will be assumed henceforth to be
proper. Improper ones seem to be of no interest.
\section{Duality of pseudocomplex structures}
We need to recall some results from McKay \cite{McKay:1999,McKay:2003}
concerning tame pseudocomplex structures on the complex
projective plane. (As always, we are assuming that our pseudocomplex
structures are proper.)
\begin{theorem}[Taubes \cite{Taubes:1995}]
There is a unique symplectic structure on the complex projective
plane, up to symplectomorphism and scaling by a constant.
\end{theorem}
\begin{theorem}[Gromov \cite{Gromov:1985}, McKay \cite{McKay:1999}]
Pick $E$ any tame pseudocomplex structure on $M=\CP{2}$. 
Define an $E$-line to be a smooth immersed sphere which is an
$E$-curve and lives in the homology class generating 
$H_2\left(\CP{2}\right).$ Every $E$-line is embedded.
Any pair of $E$-lines meet transversely at a unique point.
(Proof is a Chern class argument, similar to those below.)
Deformations of $E$-lines are unobstructed (proof: same elliptic 
theory as for the usual complex structure, following Fredholm.)
Let $M^*=M^*(E)$ be the set of $E$-lines. Then $M^*$ 
is a smooth manifold, diffeomorphic to the complex projective plane.
(We will see below the technique for proving this: affine coordinates.)
There is a unique $E$-line tangent to any given $2$-plane in $M$
which belongs to $E \subset \Gro{2}{TM}$. 
The map $E \to M^*$, given by taking an $2$-plane to the
$E$-line tangent to it, is a smooth fiber bundle. 
Write the map $E \to M$ as $\pi : E \to M.$
Consider the 4-plane field $\Theta$ on $E$ determined by the
equations $\Theta_e = \pi'(e)^{-1}e.$
Write the map $E \to M^*$ as $\pi^* : E \to M^*.$
The map
\[
E \to \Gro{2}{TM^*}
\]
given by taking each $2$-plane $e \in E$
to the $2$-plane $\left(\pi^*\right)'(e) \cdot \Theta_e \subset TM^*$
is a pseudocomplex structure on $M^*.$ An $E$-line of $M^*$
is precisely a set of points of $M^*$ consisting of just those
$E$-lines in $M$ passing through a point of $M$, i.e. points 
of $M$ are lines of $M^*$ and vice versa.
\end{theorem}
\section{Affine coordinates}
Fix a smooth tame pseudocomplex structure $E$ on $\CP{2}.$
Consider two distinct $E$-lines in $\CP{2}$, say $X$ and $Y$.
They must meet at a single point transversely. Now pick any
points $\infty_X$ in $X$ and not in $Y$, and $\infty_Y$ in $Y$
and not in $X$. Given any
point $p$ of the projective plane which is not one of these two
points, draw the line through $\infty_X$ and $p$, and find that
it strikes $Y$ at a single point $y(p) \in Y$. Similarly the line
through $\infty_X$ and $p$ strikes $Y$ at a single point $x(p) \in X.$
In this way, we smoothly map $\alpha_{XY} : \CP{2} \backslash \left\{ \infty_X,
    \infty_Y \right\} \to X \times Y.$ 
\begin{lemma}
Let $\infty$ be the line through $\infty_X$ and $\infty_Y.$
This map $\alpha_{XY}$, 
called an \emph{affine chart}, is a local diffeomorphism
\[
\alpha_{XY} : \CP{2} \backslash \infty \to X \times Y.
\]
\end{lemma}
\begin{proof}
We can differentiate this map, because the infinitesimal motions of a
line are governed by pseudoholomorphic sections of the normal bundle
(holomorphic in the osculating almost complex structure), and the
deformation theory is unobstructed by Chern class calculation. We need
to show that if we move the point $p$ infinitesimally, i.e. with a
nonzero tangent vector $v \in T_p \CP{2},$ then one of the points
$x(P),y(P)$
must move by a nonzero tangent vector.
Let $L_X$ be the line through $\infty_Y$ and $p$, and $L_Y$ be the
line
through $\infty_X$ and $p$. We form the pseudoholomorphic normal 
vector fields $A_X(v)$ on $L_X$, 
and $A_Y(v)$ on $L_Y$, determined by requiring that 
$A_X(v) = 0$ at $\infty_Y$ and $A_X(v)=v$ modulo $T_p X$ at $p$.
Existence of $A_X(v)$ comes again from Riemann--Roch,
following Duistermaat. Imagine that $A_X(v)$
vanishes at $x(p)$ and that $A_Y(v)$ vanishes at $y(p).$ 
By Chern class count, if $A_X(v)$ vanishes anywhere other than
$\infty_X$, it vanishes everywhere. But it vanishes at $x(P),$
so it must vanish everywhere, so $v$ is tangent to $L_X$.
Similarly, $v$ is tangent to $L_Y.$ So $L_X=L_Y,$ and the line
through $p$ and $\infty_X$ must contain $\infty_Y.$
\end{proof}
Similarly, we can take a \emph{dual affine chart}, defined
most simply by taking the same constuction as above, and assigning
to each line $Z$ not equal to $X$ or $Y$ its point $x(Z) \in X$ 
of intersection with $X$, and its point $y(Z) \in Y$ of intersection
with $Y$.
\begin{lemma}
This map, called a \emph{dual affine chart}, is a local diffeomorphism
\[
\hat{\alpha}_{XY} : \CP{2*} \backslash \infty^* \to X \times Y
\]
where $\infty^*$ is the set of lines not striking the point 
$X \cap Y$.
\end{lemma}
\begin{proof}
This is the dual statement to the previous lemma.
\end{proof}
Using these charts, we derive the smoothness of the 
double fibration:
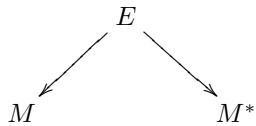
\[
\xymatrix{
 & E\ar[dl] \ar[dr] & \\
M & & M^*
}
\]
(see McKay \cite{McKay:1999})
taking a pointed $E$-line to either a point, via the left leg,
or a line via the right.
\section{Proof of the main theorem}
Take two smooth tame pseudocomplex structures $E_0$ 
and $E_1$ on $M_0=M_1=\CP{2},$
and a continuous map $\phi : M_0 \to M_1$ 
which identifies
their curves. Suppose that the map is differentiable at two
points $x_0,y_0 \in M_0$ and its derivative is an invertible
linear map.
Let $M^*_0$ and $M^*_1$ be the dual
projective planes in the two pseudocomplex structures. 
Using dual affine charts, we see that 
the map $\phi$ induces a continuous
map identifying the dual planes: call it $\phi^* : M^*_0 \to M^*_1.$ 
Then we take the lines $x^*_0, y^*_0$ in $M^*_0$
dual to the points $x_0,y_0,$ and the lines $x^*_1,y^*_1$ 
in $M^*_1$ dual to the points $x_1=\phi\left(x_0\right),
y_1=\phi\left(y_0\right).$ 
The points of the line $x^*_0$ 
are precisely the $E_0$-lines through $x_0,$ which are smoothly 
identified with their $2$-planes in $T_{x_0} M_0$, forming
the fiber of $E_0$ over $x_0.$ Under the map $\phi,$
we map $\phi'\left(x_0\right) : T_{x_0} M_0 \to T_{x_1} M_1,$ 
which is a linear map between $T_{x_0} M_0$ and $T_{x_1} M_1,$ 
so analytic. Under this analytic identification of tangent spaces,
we get an analytic identification of Grassmanians of
2-planes in those tangent spaces, and therefore a smooth
identification of the fibers of $E_0$ and $E_1$ over $x_0$
and $x_1$ (these fibers are smooth submanifolds of Grassmannians) 
so a smooth identification of the points of
the lines $x^*_0$ and $x^*_1$ Similar remarks
hold for $y^*_0$ and $y^*_1$. Finally, we use the
coordinate axis construction to produce smooth affine
coordinates on the dual projective planes, which must be
matched up by the map $\phi^*: M^*_0 \to M^*_1,$ forcing
$\phi^*$ to be smooth on the open set where the 
coordinates are defined, i.e. away from some chosen points
on the two lines. We can change the choice of those points,
and obtain global smoothness. The smoothness of the
map $\phi^*$ allows us to repeat the above argument
on the dual planes, so that $\phi$ is also smooth.  

The proof actually gives:
\begin{theorem}
A homeomorphic isomorphism of any smooth
projective planes, differentiable of full rank at two points,
is smooth.
\end{theorem}
We obtain our theorem as a consequence
of the result that smooth tame pseudocomplex
structures on $\CP{2}$ are topological projective planes, and
that the topological projective plane structure is
invariant under homeomorphic isomorphism. This
holds because the
image of an $E_0$-line is a sphere in the generating
homology class, and an $E_1$-curve, so an $E_1$-line.

%\nocite{*}
\bibliographystyle{amsplain}
\bibliography{cp2regularity}

\end{document}